\documentclass[a4paper,10pt]{article}

\usepackage{amsmath}
\def\tr{\,\mathrm{Tr\,}}

\newtheorem{theorem}{Theorem}[section]
\newtheorem{lemma}[theorem]{Lemma}

\numberwithin{equation}{section}

\begin{document}
\title 
{Dimension formula for the space of  relative symmetric polynomials of $D_n$ with respect to any irreducible representation }

\author{\bf S. Radha and P. Vanchinathan  \\
  VIT University, Chennai, India\\
\tt radha.s@vit.ac.in\quad vanchinathan.p@vit.ac.in}

\maketitle

Classification: Primary 05E05 Secondary 15A69
\date{}
\begin{abstract}
We provide an alternative formula along with the generating function for the dimension of the vector space of 
relative symmetric polynomials of $D_n$ with respect to any irreducible 
character of $D_n$.   
\end{abstract}

\section{Introduction}
Every complex representation of a finite group has a canonical decomposition
into direct sum of isotypical components.
Serre's textbook \cite[page 21]{Serre} gives the formula for the
projection map to all these components. We recall the formula here:

Given the isotypical decomposition 
$V=\bigoplus_{\chi\in \mathrm{Irr\,}(G)} V_\chi$, 
the projection to the component $V_\chi$ is given by
\begin{equation}\label{canproj}
    p= \frac{\deg \chi}{|G|} \sum_{g\in G}\bar {\chi}(g) \rho(g)
\end{equation}
When $\chi$ is the trivial 1-dimensional representation this projection
is the Reynold's operator.

In this paper we focus on the natural action of  $S_n$ (and its subgroups)
on the
complex polynomial algebra of $n$ variables, by permuting the
variables.  

We denote by $H_d(x_1,x_2,\ldots, x_n)$ the complex vector space
of all homogeneous polynomials of degree $d$ in the $n$ variables,
$x_1,x_2,\ldots, x_n$, sometimes denoted simply by $H_d$.

The image of the Reynold's operator will be the space of all symmetric
polynomials of degree $d$.

M. Shahryari \cite{Shah1} has introduced the notion of relative symmetric 
polynomials for any subgroup $G \subset S_n$ with respect to any 
irreducible character $\chi$ of $G$.

The vector space of relative symmetric polynomials  of $G$ 
relative to $\chi$, denoted by $H_d(G,\chi)$ is defined as the 
image of the  projection operator defined in Equation \eqref{canproj}  
in this case: %%$T(G,\chi)$:
\begin{equation}\label{relproj}
    T(G,\chi) = \frac{\chi(1)} {|G|} \sum_{g \in G}\chi(g) g    
\end{equation}

Finding the dimension of this vector space of relative symmetric polynomials for various subgroups of the symmetric group of $S_n$ is a fundamental question.

 In  later papers Babaei, Zamani and Shahryari found the dimension
 of the space of relative invariants for $S_n$ and its subgroup $A_n$ \cite{Shah1} and for Young subgroup \cite{Shah2}

In a series of papers Babaei and Zamani have given corresponding formula for the cyclic group in \cite{ZB1}, for the dicyclic group in \cite{ZB2} and for the dihedral group $D_n$ \cite{ZB3}. 

In our work here we relook at the formula for dihedral groups.
Babaei-Zamani formula is a summation involving cosine values which are in general irrational numbers. 
So they may be inconvenient to calculate  dimensions which are non-negative integers.
By using appropriate theorems from elementary number theory along with combinatorial counting arguements, we 
have provided an alternative dimensional formula as a summation of integer terms. 
Another advantage of our formula is that it allows us to write down the generating function.
See theorems \ref{thmdeg2} and \ref{thmdeg1} for the precise statements.

A notable feature is that this formula gives this dimension 
as a \textit{summation over the divisors of $n$}
involving the M\"obius and Euer $\phi$-functions, even though this problem
is not apparently connected with number theory.

Our paper is organized  this way: after this introduction, in the second
section a few  number-theoretic preliminaries are assembled, 
in Section 3  the character table of the dihedral group is given following  
Serre \cite{Serre}. 
In the fourth section we state and prove
the formula. Finally Section 5 illustrates how easy it is to compute  with
these formul\ae.
\section{Preliminaries}

First we set up the notations: 
\begin{itemize}
\item  We use the standard notations $\phi(n)$ and $\mu(n)$ respectively for the Euler's totient function and 
M\"obius function.
\item  For $n$ any positive integer and $r$ a divisor of $n$ we 
 denote by $S_r(n)$ the set of integers between 1 and $n$ having $r$ as their
$\gcd$. 
\[S_r(n):=\{k:1 \leq k \leq n,\gcd(k,n)=r\}\]  
\end{itemize}

We state below, without proofs, some 
 well known facts from elementary number theory as lemmas. 
These results were known to Ramanujan \cite{Ramanujan} and von Sterneck
\cite{Von}.
\begin{lemma} \label{lemma1}
With the notation as above, we have  $S_r(n)=\{rk:k \in S_1(n/r)\}$.
I.e., $S_r(n)=rS_1(n/r)$. In particular, $|S_r(n)|=|S_1(n/r)|$.  
 \end{lemma}

\begin{lemma} \label{lemma2}
The sum of all the   primitive $n^{\rm th}$ roots of unity is $\mu(n)$.
\[
 (ie) \sum_{\textstyle{k=1 \atop \gcd(k,n)=1}} ^n
    %%e^{\textstyle {\frac{2\pi ik}{n}}} 
    \exp (\textstyle {\frac{2\pi ik}{n})}
    = \mu(n)
\]
 \end{lemma}

In fact we need the following variation of Lemma 2:
\begin{lemma} \label{lemma3}
 \[ \sum_{\textstyle{k=1 \atop \gcd(k,n)=1}} ^n 
%%\cos \frac {2\pi k}n = \mu(n)
    \cos\, (2\pi k/n) = \mu(n)
\]
\end{lemma}
    \begin{lemma}\label{lemmarjsum}
        For any two positive integers $n$ and $m$
        \[\sum_{k=1}^n \exp\big({\textstyle {\frac{2\pi i mk}{n}}}\big) =\mu \Big(\frac
        n{\gcd(m,n)}\Big) \frac{\phi(n)}{\phi\big(n/\gcd(m,n)\big)}\]
    \end{lemma}
%Lemma 3 follows from Lemma 2 because: \textit{if the sum of a certain set of complex numbers is real, then 
%it coincides with  the sum of their real parts.} \\

\section{Characters of the Dihedral group $D_n$}

We write  the elements of $D_n$, as  
$D_n=\{1,\sigma,\sigma^2,\ldots,\sigma^{n-1},\tau,\tau \sigma,\tau \sigma^2,\ldots,\tau
\sigma^{n-1}\}$
  The dihedral group $D_n$  has only degree 1 and degree 2 irreducible represenations. 

\subsection{One-dimensional representations:}
\begin{itemize}
 \item When $n$ is odd there are two irreducible representations of degree 1 namely $\chi_1$ and $\chi_2$ and the 
 character table for those represenations is given below:

\begin{center} 
\begin{tabular}{|l|l|l|}
\hline 
Character & $\sigma^k$ & $\tau \sigma^k$   \\ \hline
$\chi_1$ & $1$ & $1$  \\ \hline
$\chi_2$ & $1$ &$-1$  \\ \hline 
\end{tabular}
\end{center}  
\item When $n$ is even there are four irreducible representations of degree 1 namely $\chi_1$, $\chi_2$ ,$\chi_3$ and 
$\chi_4$ and the character table for those represenations is given below:

\begin{center} 
\begin{tabular}{|l|l|l|}
\hline 
Character&$\sigma^k$ & $\tau \sigma^k$  \\ \hline
$\chi_1$ & $1$ &$1$  \\ \hline
$\chi_2$ & $1$ &$-1$ \\ \hline
$\chi_3$ & $(-1)^k$ &$(-1)^k$  \\ \hline
$\chi_4$ & $(-1)^k$ & $(-1)^{k+1}$  \\ \hline
\end{tabular} 
\end{center}
\end{itemize}

\subsection{Two-dimensional representations} 
%%Put $\zeta=e^{\frac{2\pi i}{n}}$ and
Let $h$ be a positive integer
with $h< n/2$.
 A representation $\rho_h$ of $D_n$ has the character given by
$\psi_h(\sigma^k)=2 \cos \frac{2\pi hk}{n} $ and  $\psi_h(\tau \sigma^k)=0.$

\noindent \textit{In our paper we follow the  convention that whenever
n or r is not a positive integer the binomial coefficient 
${n \choose r}$ is interpreted as zero.}
Now we can state the main result of our paper:
\section{Dimension formul\ae\  and generating functions}
\begin{theorem} \label{thmdeg2}
 Let $\psi_h$ be the irreducible  character of degree 2 of the dihedral group
    $D_n$ as above.
 Then the dimension of $H_d(D_n,\psi)$, the vector space of relative symmetric
    polynomials is described in two cases:\\

\noindent Case(i)
  h is coprime to n:
 \[
  \dim H_d(D_n,\psi_h)=\frac{2}{n}\sum_{r|n} {r+\frac{d}{n/r}-1 \choose r-1} \mu\left(\frac{n}{r}\right)   
 \] 
 The generating function in this case is given by
 \[
  \sum_{d=0}^\infty \dim H_d(D_n,\psi_h)t^d
    =\frac{2}{n}\sum_{r|n} \mu \left(\frac{n}{r}\right) (1-t^{\frac{n}{r}})^{-r}
  \]

\noindent Case(ii)
 h is not coprime to n:
 \[
    \dim H_d(D_n,\psi_h)
    =\frac{2}{n}\sum_{r|n} {r+\frac{d}{n/r}-1 \choose r-1} \mu\left(\frac{n/r}{g}\right) 
  \frac{\phi\left(\frac{n}{r}\right)} {\phi\left(\frac{n/r}{g}\right)}
 \] 
 \\
where $\mu(n)$ is the M\"obius function and $g=\gcd (h,\frac{n}{r})$.

The generating function in this case is given by
 \[
  \sum_{d=0}^\infty \dim H_d(D_n,\psi_h)t^d=\frac{2}{n}\sum_{r|n} \mu \left(\frac{n/r}{g}\right)\frac{\phi\left(\frac{n}{r}\right)}
  {\phi\left(\frac{n/r}{g}\right)} (1-t^{\frac{n}{r}})^{-r}  
  \]
 \end{theorem}

Proof:

\textit{It suffices to prove the formula for the dimension of $H_d(D_n,\psi)$ for a general $d$. The formula for the generating function is a straight forward consequence.}

Case (i):
For definiteness we fix the embedding of $D_n$ in $S_n$ with the generators of $D_n$ as below:
$D_n=\langle\sigma,\tau\rangle$ where $\sigma$ is the $n$-cycle given by (1 2 3 $\ldots n)$ and
    $\tau(j)=n+1-j$ is the  reversal permutation. 
    \textit{In fact, $D_n$ cane be embedded in $S_n$ uniquely upto conjugacy.} Now in the case of a 2-dimensional irreducible character $\chi$  of $D_n$, $\psi(\tau \sigma ^k)=0$ for all $k$. So the dimension formula reduces to the summation over the cyclic
subgroup of all rotations in $D_n$.
\begin{equation}\label{dfrot}
 \dim H_d(D_n,\psi)=\frac{\psi(1)}{|D_n|}\sum^n_{k=1}
    \psi(\sigma^k)\tr(\sigma^k)=\frac{2}{2n}\sum_{k=1}^n 2\cos\frac{2 \pi
    k}{n}\tr(\sigma^k)
\end{equation}

 Note that $\tr (\sigma^k$) is the trace of the $k^{\rm th}$ power of the 
    $n$-cycle $\sigma$ in the vector space of homogeneous polynomials
 in $n$ variables, with $\sigma$ permuting the variables cyclically. As  this vector space has all monomials of degree $d$ in $n$
 variables as basis its dimension  is ${n+d-1 \choose n-1}$.
    Being a permutation action $\tr(\sigma^k)={}$the number of monomials of
    degree $d$ in $n$ variables fixed by $\sigma^k$.
So the calculation boils down to finding the number of invariant monomials of
    degree $d$. To calculate $\tr(\sigma^k)$, let $r=\gcd(n,k)$. Then
$\sigma^k$ decomposes into a product of $r$ number of disjoint cycles of
    length $\frac{n}{r}$. For a monomial to be invariant under
$\sigma^k$, degree of all the variables within an $\frac{n}{r}$-cycle should be constant. Call these degrees $d_1,d_2,\ldots,d_r$.

\indent (ie)
\[
 d=\frac{n}{r}d_1+\frac{n}{r}d_2+\ldots+\frac{n}{r}d_r  
\]
Therefore,
\[
 d_1+d_2+\ldots+d_r=\frac{d}{n/r}
\]
A necessary condition is $d$ must be a multiple of $\frac{n}{r}$.
\textit{Let us assume this holds.}
So $\tr(\sigma^k)=$ number of ordered partitions of $\frac{d}{n/r}$ into $r$ parts. 

 This is well known to be $\displaystyle {r+\frac{d}{n/r}-1 \choose r-1}$. Substituting this value of trace in to (\ref{dfrot}), we
get 
\begin{equation}
 \dim H_d(D_n,\psi)=\frac{1}{n}\sum_{k=1}^n 2 {r+\frac{d}{n/r}-1 \choose r-1} \cos\frac{2\pi k}{n} 
\end{equation}

Note that for  two terms of the summation on 
righthand side  if $1 \leq  k_1,k_2 \leq n$ are such that
$\gcd(k_1,n)=\gcd(k_2,n)$, 
then the coefficient of $\cos\frac{2\pi k_1}{n}$ equals that of
$\cos\frac{2\pi k_2}{n}$. As the $\gcd$ of 
any number with $n$ is a divisor of $n$, the dimension formula can be rewritten as a summation over the
divisors of $n$.
We treat the above summation as a sum of binomial coefficients ${r+\frac{d}{n/r}-1 \choose r-1}$, one for each divisor $r$ of $n$ with some weights.
These weights are sums of cosine values. So the dimensional formula takes the form 
\begin{equation}
 \dim H_d(D_n,\psi)=\frac{1}{n}\sum_{r|n} \Big[\sum_{k \in S_r(n)} 
    2 \cos\, \big( 2\pi k/n\big)
    %2 \cos\frac{ 2\pi k}{n}
    \Big] {r+\frac{d}{n/r}-1 \choose r-1}
\end{equation}
That is,
\begin{equation}
 \dim H_d(D_n,\psi)=\frac{1}{n}\sum_{r|n} \Big[\sum_{a \in S_1(n/r)} 
    2 \cos\, \big( 2\pi a \big/{\textstyle\frac nr}\big)
    %2 \cos\frac{ 2\pi a}{n/r}
    \Big] {r+\frac{d}{n/r}-1 \choose r-1} \mbox{(by lemma \ref{lemma1})}
\end{equation}

Using lemma \ref{lemma3} the above equation reduces to
\[
\dim H_d(D_n,\psi)=\frac{2}{n} \sum_{r|n} \mu\left(\frac{n}{r}\right) {r+\frac{d}{n/r}-1 \choose r-1} 
\] 

Proof of case (ii):
Proceeding as in case (i), we have
\begin{equation} \label{hdformula}
  \dim H_d(D_n,\psi)=\frac{2}{n}\sum_{r|n} \left[\sum_{a \in S_1(n/r)}  \cos\frac{ 2\pi ah}{n/r}\right] {r+\frac{d}{n/r}-1 \choose r-1} 
 \end{equation}
 The inner summation inside the square brackets in the above equation  is actually the sum of the real parts of $h^{\rm th}$ powers of all primitive $\left(\frac{n}{r}\right)^{\rm th}$
 roots of unity. Defining $g=\gcd(h,n/r)$, we see that the above is same as the sum of all the real parts of $g^{\rm th}$ powers of 
 all primitive $\big(\frac{n}{r}\big)^{\rm th}$  roots of unity. 
 Now we can apply Lemma  \ref{lemmarjsum} which makes
  equation \eqref{hdformula} to become
 \[
  \dim H_d(D_n,\psi)=\frac{2}{n}\sum_{r|n} \left[\sum_{a \in S_1(n/gr)}  \cos\frac{ 2\pi a}{n/gr}\right] \frac{\phi(n/r)}{\phi(n/gr)} {r+\frac{d}{n/r}-1 \choose r-1} 
 \]
 Again using lemma \ref{lemma3}
 \[
   \dim H_d(D_n,\psi)=\frac{2}{n}\sum_{r|n} \mu\left(\frac{n/r}{g}\right) \frac{\phi(n/r)}{\phi(n/gr)}{r+\frac{d}{n/r}-1 \choose r-1} 
 \]
\begin{theorem} \label{thmdeg1}
Case (i),  $n$ is odd: \\
Let $\chi_1$ and $\chi_2$ be the two irreducible characters of degree 1 with
    $\chi_1$ being the trivial character and $\psi_2$ taking $+1$ on rotations and $-1$ on reflections.
The dimensions of $H_d(D_n,\chi_1)$ and $H_d(D_n,\chi_2)$ are given by
 \[
  \dim H_d(D_n,\chi_1)=\frac{1}{2n}\left[ \sum_{r|n} {r+\frac{d}{n/r}-1 \choose
    r-1}\phi\left(\frac{n}{r}\right) +n \sum_{l=0}^{ \lfloor d/2 \rfloor} {\frac{n-1}{2}+l-1 \choose l}\right]
 \]
 and

\[
  \dim H_d(D_n,\chi_2)=\frac{1}{2n}\left[ \sum_{r|n} {r+\frac{d}{n|r}-1 \choose r-1}\phi\left(\frac{n}{r}\right) -n \sum_{l=0}^{\lfloor d/2 \rfloor} {\frac{n-1}{2}+l-1 \choose l}\right]
\]

The generating functions for the above two cases are given by
\[
  \sum_{d=0}^\infty \dim H_d(D_n,\chi_1)t^d=\frac{1}{2n}\left[\sum_{r|n} \phi \left(\frac{n}{r}\right)(1-t^{\frac{n}{r}})^{-r} +\frac{ n(1-t^2)^{-(n-1)/2} } {1-t}\right]
 \]
\[
  \sum_{d=0}^\infty \dim H_d(D_n,\chi_2)t^d=\frac{1}{2n}\left[\sum_{r|n} \phi \left(\frac{n}{r}\right)(1-t^{\frac{n}{r}})^{-r} -\frac{ n(1-t^2)^{-(n-1)/2} } {1-t}\right]
 \]
Case (ii) when $n$ is even: \\
Let $\chi_1$,$\chi_2$,$\chi_3$ and $\chi_4$ be the four irreducible characters of degree 1. Let $H_d(D_n,\psi_1),$   $ H_d(D_n,\psi_2), H_d(D_n,\psi_3)$  and
$H_d(D_n,\chi_4)$ be the space of Relative symmetric polynomials with respect to $\psi_1$,$\psi_2$,$\psi_3$ and $\psi_4$. Then the dimensions 
$H_d(D_n,\chi_1),$ $H_d(D_n,\chi_2),\allowbreak H_d(D_n,\chi_3)$ and $H_d(D_n,\chi_4)$ are given by
\begin{equation*}
\begin{aligned}
  \dim H_d(D_n,\chi_1)&=\frac{1}{2n}\left\{\sum_{r|n}{r+\frac{d}{n/r}-1 \choose r-1}\phi\left(\frac{n}{r}\right) \right.\\
                       & {}\qquad \left.+\frac{n}{2}\left[{\frac{n}{2}+\frac{d}{2}-1 \choose \frac{d}{2}}+\sum_{l=0}^{ d/2}{\frac{n-2}{2}+l-1 \choose l}(d-2l+1)\right]\right\}
\end{aligned}
\end{equation*}
\begin{equation*}
\begin{aligned}
\dim H_d(D_n,\chi_2)&=\frac{1}{2n}\left\{\sum_{r|n}{r+\frac{d}{n/r}-1 \choose r-1}\phi\left(\frac{n}{r}\right) \right. \\
                     & {}\qquad \left. -\frac{n}{2}\left[{\frac{n}{2}+\frac{d}{2}-1 \choose \frac{d}{2}}+\sum_{l=0}^{d/2}{\frac{n-2}{2}+l-1 \choose l}(d-2l+1)\right]\right\}
\end{aligned}
\end{equation*}
\begin{equation*}
\begin{aligned}
\dim H_d(D_n,\chi_3) &= \frac{1}{2n}\left \{\sum_{r|n} (-1)^r \phi\left(\frac{n}{r}\right){r+\frac{d}{n/r}-1 \choose r-1} \right. \\
                     &{}\qquad+\left.\frac{n}{2}\left [{\frac{n}{2}+\frac{d}{2}-1 \choose \frac{d}{2}} - \sum_{l=0}^{ \lfloor d/2 \rfloor}{\frac{n-2}{2}+l-1 \choose l}(d-2l+1) \right] \right\}
\end{aligned}
\end{equation*}
\begin{equation*}
\begin{aligned}
\dim H_d(D_n,\chi_4) & = \frac{1}{2n}\left \{\sum_{r|n} (-1)^r \phi\left(\frac{n}{r}\right){r+\frac{d}{n/r}-1 \choose r-1} \right. \\
                     & {}\qquad - \left.\frac{n}{2}\left [{\frac{n}{2}+\frac{d}{2}-1 \choose \frac{d}{2}} - \sum_{l=0}^{\lfloor d/2 \rfloor}{\frac{n-2}{2}+l-1 \choose l}(d-2l+1) \right] \right\}
\end{aligned}
\end{equation*}
The generating functions for the above four cases are given by:
\[
  \sum_{d=0}^\infty \dim H_d(D_n,\chi_1)t^d=\frac{1}{2n}\bigg[\sum_{r|n} \phi
    \left(\frac{n}{r}\right)(1-t^{\frac{n}{r}})^{-r} +\frac{n}{2}(1-t^2)^{-(n+2)/2}
    (2+t^2)(1+t)\bigg]
 \]
\[
  \sum_{d=0}^\infty \dim H_d(D_n,\chi_2)t^d=\frac{1}{2n}\left[\sum_{r|n} \phi \left(\frac{n}{r}\right)(1-t^{\frac{n}{r}})^{-r}-\frac{n}{2}(1-t^2)^{-(n+2)/2} (2+t^2)(1+t)\right]
 \]
\begin{equation*}
\begin{aligned}
\sum_{d=0}^\infty & \dim H_d(D_n,\chi_3)t^d \\
                  &{}\qquad =\frac{1}{2n}\left[\sum_{r|n} \phi \left(\frac{n}{r}\right)(-1)^r(1-t^{\frac{n}{r}})^{-r} - \frac{n}{2}(1-t^2)^{-(n+2)/2} (1+t+t^2)\right]
\end{aligned}
 \end{equation*}
\begin{equation*}
\begin{aligned}
\sum_{d=0}^\infty &\dim H_d(D_n,\chi_4)t^d \\
&{}\qquad =\frac{1}{2n}\left[\sum_{r|n} \phi \left(\frac{n}{r}\right)(-1)^r(1-t^{\frac{n}{r}})^{-r} + \frac{n}{2}(1-t^2)^{-(n+2)/2} (1+t+t^2)\right]
\end{aligned}
\end{equation*}
 \end{theorem}

Proof:

%\textit{In the same way it suffices to prove the formula for the dimension of
%$H_d(D_n,\chi)$ for a general $d$. The formula for the generating function is a
%simple consequence.}\\

Case (i) (when $n$ is odd): \\
By definition:
\[
 \dim H_d(D_n,\chi_1)= \frac{1}{2n}\left[\sum_{k=1}^n \tr(\sigma ^k) \psi_1(\sigma ^k) + \sum_{k=1}^n \tr(\tau \sigma ^k) \psi_1(\tau \sigma ^k)\right]
\]
(ie)
\[
 \dim H_d(D_n,\chi_1)= \frac{1}{2n}\left[\sum_{k=1}^n \tr(\sigma ^k)  + \sum_{k=1}^n \tr(\tau \sigma ^k)\right]
\]
When $n$ is odd all the reflections of $D_n$ falls into a single conjugacy class. Hence the above summation becomes
\[
 \dim H_d(D_n,\chi_1)= \frac{1}{2n}\left[n \tr(\tau) + \sum_{k=1}^n \tr(\sigma ^k) \right]
\]
Now calculation of $\tr(\sigma^k)$ is the same as in Theorem 1. It remains to find $\tr(\tau)$. 
Since $n$ is odd, $\tau$ is the product of $\frac{n-1}{2}$ transpositions and
has one fixed point.
Now $\tr(\tau)$ is the count of monomials fixed by $\tau$. We denote the variables by $x_1,x_2,\ldots x_{(n-1)/2},y_1,y_2,\ldots,y_{(n-1)/2},z$. 
Without loss of generality, let us assume that $\tau(z)=z$, $\tau(x_i)=y_i$ and $\tau(y_i)=x_i$ for $1 \leq i \leq (n-1)/2$.
Now a monomial of degree $d$ invariant under $\tau$ has to be of the form
\[
 z^{d_0} x_1^{d_1} y_1^{d_1} x_2^{d_2} y_2^{d_2}  x_3^{d_3} y_3^{d_3}\ldots 
\]

Number of tuples $(d_0,d_1,\ldots,d_{(n-1)/2})$ such that
\[
 d_0+2(d_1+d_2+\ldots+d_{(n-1)/2})=d
\]
is the total number of ordered partitions of $(d-d_0)/2$ into $(n-1)/2$ parts with all $d_0$ satisfying $d-d_0$ is an even non negative integer.
This is easily verified to be $\sum_{l=0}^{ \lfloor d/2 \rfloor} {\frac{n-1}{2}+l-1 \choose l}$. Hence the formula becomes

\[
 \dim H_d(D_n,\chi_1)= \frac{1}{2n}\left[\sum_{r|n} {r+\frac{d}{n/r}-1 \choose r-1} +n \sum_{l=0}^{\lfloor d/2 \rfloor} {\frac{n-1}{2}+l-1 \choose l}\right]
\]
Using the same arguements as above, we have
\[
 \dim H_d(D_n,\chi_2)= \frac{1}{2n}\left[\sum_{r |n} {r+\frac{d}{n/r}-1 \choose r-1} -n \sum_{l=0}^{ \lfloor d/2 \rfloor} {\frac{n-1}{2}+l-1 \choose l}\right]
\]
case (ii):
By definition 
\[
 \dim H_d(D_n,\chi_1)= \frac{1}{2n}\left[\sum_{k=1}^n \tr(\sigma ^k) \psi_1(\sigma ^k) + \sum_{k=1}^n \tr(\tau \sigma ^k) \psi_1(\tau \sigma ^k)\right]
\]
(ie)
\[
 \dim H_d(D_n,\chi_1)= \frac{1}{2n}\left[\sum_{k=1}^n \tr(\sigma ^k)  + \sum_{k=1}^n \tr(\tau \sigma ^k)\right]
\]
Now $\tr(\sigma^k)$ is the same as in Theorem 1. It remains to find $\tr(\tau \sigma^k)$. Since $n$ is even, all the reflections 
$\tau \sigma^k$ fall into two conjugacy classes according as  $k$ is even or odd. Hence the formula becomes
\[
 \dim H_d(D_n,\chi_1)= \frac{1}{2n}\left[\sum_{k=1}^n \tr(\sigma ^k)  + \frac{n}{2}\tr(\tau \sigma)+\frac{n}{2} \tr(\tau )\right]
\]
Using the same arguement as we did in case (i) of this theorem, we can find the trace of the reflections in both the 
cases (when $n$ is odd or even). Hence the formula becomes

\begin{equation*}
\begin{aligned}
\dim H_d(D_n,\chi_1)&=\frac{1}{2n}\left\{\sum_{r|n}{r+\frac{d}{n/r}-1 \choose r-1} \phi\left(\frac{n}{r}\right) \right. \\
& {}\qquad \left. +\frac{n}{2}\sum_{l=0}^{d/2} \left[{\frac{n-2}{2}+l-1 \choose l}(d-2l+1)+{\frac{n}{2}+l-1 \choose l} \right]     \right\}
\end{aligned}
\end{equation*}
In the same way 
\begin{equation*}
\begin{aligned}
\dim H_d(D_n,\chi_2)& =\frac{1}{2n}\left\{\sum_{r|n}{r+\frac{d}{n/r}-1 \choose r-1} \phi\left(\frac{n}{r}\right) \right.\\
&{} \qquad \left. -\frac{n}{2}\sum_{l=0}^{\lfloor d/2 \rfloor}\left[{\frac{n-2}{2}+l-1 \choose l}(d-2l+1)+{\frac{n}{2}+l-1 \choose l} \right] \right\}
\end{aligned}
\end{equation*}

Now again using the definition, we have
\[
 \dim H_d(D_n,\chi_3)= \frac{1}{2n}\left[\sum_{k=1}^n \tr(\sigma ^k) \psi_3(\sigma ^k) + \sum_{k=1}^n \tr(\tau \sigma ^k) \psi_3(\tau \sigma ^k)\right]
\]
(ie)\[
 \dim H_d(D_n,\chi_3)= \frac{1}{2n}\left[\sum_{k=1}^n (-1)^k) \tr(\sigma ^k) + \sum_{k=1}^n (-1)^k \tr(\tau \sigma ^k) \right]
\]
Using the same result for $\tr(\sigma^k)$ from Theorem 1, we have
\[
 \dim H_d(D_n,\chi_3)= \frac{1}{2n}\left\{ \sum_{r|n} \left[ \left(\sum_{k \in S_r(n)} (-1)^k \right) {r+\frac{d}{n/r}-1 \choose r-1} \right] + \sum_{k=1}^n (-1)^k \tr(\tau \sigma ^k) \right\}
\]
Since $n$ is even, all $k$'s in the inner summation are even or odd according as $r$ is even or odd and it sums up to $(-1)^r \phi(\frac{n}{r})$.
Hence the above equation becomes
\[
 \dim H_d(D_n,\chi_3)= \frac{1}{2n}\left\{ \sum_{r|n} \left[ (-1)^r \phi\left(\frac{n}{r}\right){r+\frac{d}{n/r}-1 \choose r-1} \right] + \sum_{k=1}^n (-1)^k \tr(\tau \sigma ^k) \right\}
\]
It remains to find $\tr(\tau \sigma ^k)$. Using the same arguement as we did in case (i), we can calculate it. Hence the above equation reduces to
\begin{equation*}
\begin{aligned}
\dim H_d(D_n,\chi_3) & = \frac{1}{2n} \left\{\sum_{r|n}(-1)^r \phi \left(\frac{n}{r}\right){r+\frac{d}{n/r}-1 \choose r-1} \right. \\
                     &{}\qquad \left. +\frac{n}{2} \left[{\frac{n}{2}+\frac{d}{2}-1 \choose \frac{d}{2}}-\sum_{l=0}^{d/2}{\frac{n-2}{2}+l-1 \choose l}(d-2l+1)\right] \right\}
\end{aligned}
\end{equation*}
In the same way, we have
\begin{equation*}
\begin{aligned}
 \dim H_d(D_n,\chi_4) &= \frac{1}{2n}\left \{\sum_{r|n} (-1)^r \phi\left(\frac{n}{r}\right){r+\frac{d}{n/r}-1 \choose r-1} \right. \\
                      & \left.-\frac{n}{2}\left [{\frac{n}{2}+\frac{d}{2}-1 \choose l} - \sum_{l=0}^{\lfloor d/2 \rfloor}{\frac{n-2}{2}+l-1 \choose l}(d-2l+1) \right] \right\}
\end{aligned}
\end{equation*}

\textit{To find the generating functions for $\dim H_d(D_n,\chi_1),$ $\dim H_d(D_n,\chi_2)$, $\dim H_d(D_n,\chi_3)$ and $\dim H_d(D_n,\chi_4)$:}

Let $G_1(t)$ be the generating function for $\sum_{r|n} \phi \left(\frac{n}{r}\right){r+ \frac{d}{n/r}-1 \choose r-1}$ which is easily seen to be 
$\sum_{r|n} (-1)^r \phi \left(\frac{n}{r}\right)(1-t^{\frac{n}{r}})^{-r}$. Let $G_2(t)$ and $G_3(t)$ be the generating functions corresponding to 
${\frac{n}{2}+\frac{d}{2}-1 \choose l}$ and $\sum_{l=0}^{d/2}{\frac{n-2}{2}+l-1 \choose l}(d-2l+1)$ respectively. Now
$G_2(t)$ is very easily verified to be $(1-t^2)^{\frac{n}{2}}$. 
Now the generating functions for $\dim(H_d,\chi_1)$ ,$\dim(H_d,\chi_2)$ ,$\dim(H_d,\chi_3)$ and $\dim(H_d,\chi_4)$ are given by:
\[
\sum_{d=0}^{\infty} \dim H_d(D_n,\chi_1)t^d=\frac{1}{2n}\left\{G_1(t)+\frac{n}{2}\left[G_2(t)+G_3(t)\right]\right\}
\]
\[
\sum_{d=0}^{\infty} \dim H_d(D_n,\chi_2)t^d=\frac{1}{2n}\left\{G_1(t)-\frac{n}{2}\left[G_2(t)+G_3(t)\right]\right\}
\]
\[
\sum_{d=0}^{\infty} \dim H_d(D_n,\chi_3)t^d=\frac{1}{2n}\left\{G_1(t)+\frac{n}{2}\left[G_2(t)-G_3(t)\right]\right\}
\]
\[
\sum_{d=0}^{\infty} \dim H_d(D_n,\chi_3)t^d=\frac{1}{2n}\left\{G_1(t)-\frac{n}{2}\left[G_2(t)-G_3(t)\right]\right\}
\]
Now evaluation of $G_3(t)$ given in the following lemma.

\begin{lemma}
For 
\[
\sum_{l=0}^{\lfloor d/2 \rfloor}{\frac{n-2}{2}+l-1 \choose l}(d-2l+1)
\]
the generating function is
\[
G_3(t)=(1-t^2)^{-\frac{n+2}{2}}(1+2t+t^2+t^3)
\]
\end{lemma}
Proof of Lemma: 

First consider
\[
 \sum_{l=0}^{\lfloor d/2 \rfloor}{\frac{n-2}{2}+l-1 \choose l}(d-2l+1)
\]
Let 
\[
 b_d=\sum_{l=0}^{\lfloor d/2 \rfloor}{m+l-1 \choose l}(d-2l+1)
\]
One easily verifies 

\[
 b_{2d+1}=b_{2d}+\sum_{l=0}^d {m+1-1 \choose l}
\]
Now again induction gives

\[
 \sum_{l=0}^d {m+1-1 \choose l}={m+d \choose d}  
\]
Therefore,

\[
 b_{2d+1}=b_{2d}+ {m+d \choose d}
\]
Let
\[
\sum_{d=0}^{\infty} b_{2d} t^{2d}=G_5(t)
\]
Now

\[
\sum_{t=0}^{\infty}b_{2d+1}t^{2d+1}=t\sum_{t=0}^{\infty} \left[b_{2d} t^{2d}+{m+d \choose d} t^{2d} \right]=tG_5(t)+t(1-t^2)^{-(m+1)}
\]
Hence 

\[
    G_4(t)=\sum_{t=0}^{\infty}b_dt^d=\sum_{t\, \mathrm{even}}+\sum_{t,\mathrm{odd}}=G_5(t)+tG_5(t)+t(1-t^2)^{-(m+1)}
\]
To evaluate $G_5(t)$

\begin{equation*}
\begin{aligned}
 G_5(t)&=\sum_{d=0}^{\infty} \sum_{l=0}^d {m+l-1 \choose l} (2d-2l+1)t^{2d} \\
      &=\sum_{d=0}^{\infty} \sum_{l=0}^d {m+l-1 \choose l} (2d+1)t^{2d}-2 \sum_{l=0}^d l {m+l-1 \choose l}t^{2d}
\end{aligned}
\end{equation*}

Now using induction one can show that
\begin{equation*}
\begin{aligned}
 \sum_{l=0}^d l {m+l-1 \choose l}&=m \sum_{l=0}^d \frac{l+m-m}{m} {m+l-1 \choose l} \\
&=m\sum_{l=0}^d \frac{m+l}{m} {m+l-1 \choose l} - m \sum_{l=0}^d{m+l-1 \choose l}
\end{aligned}
\end{equation*}
(ie)
\begin{equation*}
\begin{aligned}
 \sum_{l=0}^d l {m+l-1 \choose l}&=m \sum_{l=0}^d\frac{m+l}{m}\frac{(m+l-1)!}{l! (m-1)!}-m\sum_{l=0}^d{m+l-1 \choose l} \\
&=m\sum_{l=0}^d{m+l \choose l}- m \sum_{l=0}^d{m+l-1 \choose l}
\end{aligned}
\end{equation*}
After simplification,we have
\[
\sum_{l=0}^d l {m+l-1 \choose l}=m{m+d \choose d-1}
\]
Hence,
\begin{equation*}
\begin{aligned}
 G_5(t)&=\sum_{d=0}^{\infty} {m+d \choose d} (2d+1)t^{2d}-2 \sum_{d=0}^{\infty} m {m+d \choose d-1}t^{2d} \\
       &=\sum_{d=0}^{\infty} {m+d \choose d} 2dt^{2d}+\sum_{d=0}^{\infty} {m+d \choose d}t^{2d}-2 \sum_{d=0}^{\infty} m{m+d \choose d-1}t^{2d}
\end{aligned}
\end{equation*}

Now using
\[
\sum_{d=0}^{\infty} {m+d \choose d} 2d t^{2d}=t \frac{d}{dt}\left[(1-t^2)^{-(m+1)}\right]
\]
and 
\[
 \sum_{d=0}^{\infty} m{m+d \choose d}t^{2d}=2mt^2(1-t^2)^{-(m+2)}
\]
Simplifying 
\[
 G_5(t)=(1-t^2)^{-(m+2)}(2t^2+1)
\]

Finally 
\[
 G_4(t)=(1-t^2)^{-(m+2)}(1+2t+2t^2+t^3)
\]
Replacing $m$ by $\frac{n-2}{2}$ in $G_4(t)$, we have
\[
 G_3(t)=(1-t^2)^{-(n+2)/2}(1+2t+2t^2+t^3)
\]

\section{Examples}
\subsection{The dihedral group $D_{10}$}
Consider the group $D_{10}$. It has 4 one-dimensional irreducible representations and 4 two-dimensional irreducible representations. Here we give the generating functions for all the eight irreducible representations.

Two-dimensional representations:
As $0 < h < \frac{n}{2}$, h can assume 1,2,3 and 4

case (i) when $h$ is coprime with $n$ i.e.,  $ h=1,3$.

The generating function is
\[
\frac{1}{5}\left[\frac{1}{(1-t^{10})}-\frac{1}{(1-t^5)^2}-\frac{1}{(1-t^2)^5}+\frac{1}{(1-t)^{10}}\right]
\]

case (i) when $h$ is not coprime with $n$  i.e., $ h=2,4$

The generating function is
\[
\frac{1}{5}\left[-\frac{1}{(1-t^{10})}-\frac{1}{(1-t^5)^2}+\frac{1}{(1-t^2)^5}+\frac{1}{(1-t)^{10}}\right]
\]

One-dimensional representations:

The generating functions for $\chi_1$ and $\chi_2$ are given by
\[
\frac{1}{20}\left\{\left[\frac{1}{(1-t^{10})}+\frac{4}{(1-t^5)^2}+\frac{1}{(1-t^2)^5}+\frac{1}{(1-t)^{10}}\right] \pm 10 \left[\frac{(2+t^2)(1+t)}{(1-t^2)^{6}}\right]\right\}
\]
The generating functions for $\chi_3$ and $\chi_4$ are given by
\[
\frac{1}{20}\left\{\left[-\frac{1}{(1-t^{10})}+\frac{4}{(1-t^5)^2}-\frac{1}{(1-t^2)^5}+\frac{1}{(1-t)^{10}}\right] \mp 10 \left[\frac{(1+t+t^2)}{(1-t^2)^{6}}\right]\right\}
\]
\section{Existence of relative invariants}
The above formul\ae\  give that the dimension of the space of
relative symmetric polynomials of degree 1 is 2 whatever be the character.
In particular relative invariants for $D_n$ exist in degree 1 always.

One can easily  see that the dimension is positive for $D_n$ for any degree $d$
when $n$ a prime
number. It seems all the dimensions are always positive though we are unable to
prove this.

\end{document}